\documentclass[12pt]{article}

\usepackage[T1]{fontenc}
\usepackage[utf8]{inputenc}
\usepackage{lmodern}
\usepackage[british]{babel}
\usepackage{mathtools,amssymb,amsthm}
\usepackage[hidelinks]{hyperref}
\usepackage{geometry}
\theoremstyle{plain}
\newtheorem{theorem}{Theorem}
\newtheorem{proposition}[theorem]{Proposition}

\theoremstyle{remark}
\newtheorem{remark}[theorem]{Remark}

\newcommand{\ip}[2]{\langle #1,#2\rangle}

\title{An Identity of Fillipi Related to Fermat-Type Equations}

\author{Mike Winkler$^1$\quad Andreas Fillipi$^2$}

\date{\small
$^1$Ruhr University Bochum, Faculty of Mathematics, Germany\\
\texttt{mike.winkler@ruhr-uni-bochum.de}\\[2mm]
$^2$Achern, Germany\\
\texttt{andreas.fillipi@googlemail.com}\\[4mm]
\normalsize August 21, 2026
}

\begin{document}

\maketitle

\begin{abstract}
We give a structural proof of an algebraic identity proposed by Fillipi in 2014. For every integer $n\ge2$, the identity represents an explicit multiple of $x^n+y^n-z^n$ as $\mathcal A^2+\mathcal B^2-\mathcal C^2$. We identify this expression with a determinant arising from a weighted Gram matrix and show that Fillipi's particular fourth-power terms are obtained from a simpler identity by one polynomial column transformation. The same transformation yields an explicit polynomial rank-one factorisation modulo $x^n+y^n-z^n$. Thus the associated Pythagorean relation is a universal polynomial consequence of the Fermat-type equation and supplies no additional algebraic solvability condition.
\end{abstract}

\section{Introduction}

In 2014, Fillipi proposed an algebraic identity in connection with Fermat-type equations~\cite{Fillipi2014}. For an integer $n\ge 2$, put $k=n-2$ and
\[
r=x-y,\qquad s=y+z,\qquad t=z+x,
\]
\[
u=x+y+z,\qquad v=y-z-x,\qquad w=x-y-z.
\]
Define
\begin{align*}
P&=r^2u^4(xy)^k-s^2v^4(yz)^k-t^2w^4(zx)^k,\\
Q&=r^2(xy)^k-s^2(yz)^k-t^2(zx)^k,\\
R&=(ru)^2(xy)^k-(sv)^2(yz)^k-(tw)^2(zx)^k,\\
\mathcal A&=P-Q,\qquad
\mathcal B=2R,\qquad
\mathcal C=P+Q.
\end{align*}
The identity is
\begin{equation}\label{eq:fillipi-intro}
\mathcal A^2+\mathcal B^2-\mathcal C^2
=
(8rst)^2(xyz)^{n-2}(x^n+y^n-z^n).
\end{equation}

The determinant mechanism used below is classical. The point of this note is to identify the specific structure behind \eqref{eq:fillipi-intro} and to determine its exact algebraic content on $x^n+y^n=z^n$. We first derive a simpler identity, with factor $4r^2$, from Jacobi's complementary-minor identity for a weighted Gram matrix. A single polynomial column transformation then produces both Fillipi's factor $64r^2s^2t^2$ and the specific forms of $P$ and $R$. Finally, the same transformation gives explicit polynomial parameters for the resulting Pythagorean relation. In particular, the latter is formally forced by the Fermat-type equation and is not an additional algebraic solvability criterion.

\section{A weighted Gram determinant identity}

Let $S$ be a commutative ring with unit. Let $D=\operatorname{diag}(d_1,d_2,d_3)$ with entries in $S$, and write
\[
\ip{a}{b}_D=a^{\mathsf T}Db
\]
for column vectors $a,b\in S^3$.

\begin{proposition}\label{prop:gram}
Let $\xi,a,b\in S^3$, let
\[
F_0=\ip{\xi}{\xi}_D,
\]
and define
\[
\Gamma_D(\xi;a,b)
=
\begin{pmatrix}
P_0&R_0\\
R_0&Q_0
\end{pmatrix},
\]
where
\begin{align*}
P_0&=F_0\ip{a}{a}_D-\ip{\xi}{a}_D^2,\\
Q_0&=F_0\ip{b}{b}_D-\ip{\xi}{b}_D^2,\\
R_0&=F_0\ip{a}{b}_D-\ip{\xi}{a}_D\ip{\xi}{b}_D.
\end{align*}
Then
\begin{equation}\label{eq:gramdet}
\det\Gamma_D(\xi;a,b)
=
F_0\det(D)[\xi,a,b]^2,
\end{equation}
where $[\xi,a,b]$ denotes the determinant of the matrix with columns $\xi,a,b$.
Consequently,
\begin{equation}\label{eq:lorentz}
(P_0-Q_0)^2+(2R_0)^2-(P_0+Q_0)^2
=
-4F_0\det(D)[\xi,a,b]^2.
\end{equation}
\end{proposition}

\begin{proof}
Let $M=[\xi,a,b]$ and $G=M^{\mathsf T}DM$. The principal $2\times2$ submatrix of $\operatorname{adj}G$ complementary to the first row and column is
\[
\begin{pmatrix}
Q_0&-R_0\\
-R_0&P_0
\end{pmatrix}.
\]
Jacobi's complementary-minor identity gives
\[
P_0Q_0-R_0^2=G_{11}\det G.
\]
This is a polynomial identity in the entries of $G$ and hence remains valid
over an arbitrary commutative ring. Since $G_{11}=F_0$ and
\[
\det G=\det(D)(\det M)^2,
\]
equation \eqref{eq:gramdet} follows. Equation \eqref{eq:lorentz} is the elementary identity
\[
(P_0-Q_0)^2+(2R_0)^2-(P_0+Q_0)^2
=-4(P_0Q_0-R_0^2).
\]
\end{proof}

\begin{remark}
Proposition~\ref{prop:gram} is classical. It is Jacobi's identity for complementary minors applied to a Gram matrix, together with the determinant formula for $M^{\mathsf T}DM$; see, for example, \cite[Ch.~I, Sec.~4]{Gantmacher} and \cite[Sec.~0.8.5, Eq.~(0.8.5.8)]{HornJohnson}.
\end{remark}

For a $3\times2$ matrix $B=[a,b]$, we shall write
\[
\Gamma_D(\xi;B)=\Gamma_D(\xi;a,b).
\]
From the definition,
\begin{equation}\label{eq:matrixform}
\Gamma_D(\xi;B)
=
F_0B^{\mathsf T}DB
-
(B^{\mathsf T}D\xi)(B^{\mathsf T}D\xi)^{\mathsf T}.
\end{equation}

\section{The base identity and the Fillipi transformation}

From now on work in $\mathbb Z[x,y,z]$ and put
\[
F=x^n+y^n-z^n,
\qquad
D=\operatorname{diag}(x^k,y^k,-z^k),
\qquad
\xi=(x,y,-z)^{\mathsf T}.
\]
Further, let
\[
e=(1,1,1)^{\mathsf T},
\qquad
e_3=(0,0,1)^{\mathsf T}.
\]
Then $\ip{\xi}{\xi}_D=F$ and $\det D=-(xyz)^k$.

For the pair $(e_3,e)$, write
\[
\Gamma_0
=
\Gamma_D(\xi;e_3,e)
=
\begin{pmatrix}
P_0&R_0\\
R_0&Q
\end{pmatrix}.
\]
A direct calculation gives
\begin{equation}\label{eq:baseentries}
P_0=-z^k(x^n+y^n),
\qquad
R_0=-z^k\bigl(tx^{n-1}+sy^{n-1}\bigr),
\end{equation}
and
\begin{equation}\label{eq:Qexplicit}
Q=r^2(xy)^k-s^2(yz)^k-t^2(zx)^k.
\end{equation}
Moreover,
\[
[\xi,e_3,e]=y-x=-r.
\]
Proposition~\ref{prop:gram} therefore yields the base identity
\begin{equation}\label{eq:baseidentity}
(P_0-Q)^2+(2R_0)^2-(P_0+Q)^2
=
4r^2(xyz)^kF.
\end{equation}

Now set
\begin{equation}\label{eq:g}
g=u^2e-4st\,e_3,
\end{equation}
and define
\[
B_0=[e_3,e],
\qquad
B=[g,e].
\]
Then
\begin{equation}\label{eq:basis}
B=B_0T,
\qquad
T=
\begin{pmatrix}
-4st&0\\
u^2&1
\end{pmatrix}.
\end{equation}
This is a polynomial column transformation rather than a basis change over
$\mathbb Z[x,y,z]$, since $\det T=-4st$ need not be a unit. More precisely,
\[
T=
\begin{pmatrix}
-4st&0\\
0&1
\end{pmatrix}
\begin{pmatrix}
1&0\\
u^2&1
\end{pmatrix},
\]
so the transformation consists of a scaling of the first column followed by a
unimodular shear. By \eqref{eq:matrixform},
\begin{equation}\label{eq:congruence}
\Gamma_D(\xi;B)
=
T^{\mathsf T}\Gamma_0T.
\end{equation}
Since the second column of $T$ is $(0,1)^{\mathsf T}$, the lower-right entry is unchanged and remains $Q$.
Hence, writing
\[
\Gamma_D(\xi;B)
=
\begin{pmatrix}
P&R\\
R&Q
\end{pmatrix},
\]
we obtain
\begin{equation}\label{eq:PRtransform}
R=-4stR_0+u^2Q,
\qquad
P=16s^2t^2P_0-8stu^2R_0+u^4Q.
\end{equation}

The first formula in \eqref{eq:PRtransform}, together with \eqref{eq:baseentries} and \eqref{eq:Qexplicit}, gives
\begin{align*}
R
&=r^2u^2(xy)^k
-s^2(u^2-4ty)(yz)^k
-t^2(u^2-4sx)(zx)^k\\
&=(ru)^2(xy)^k-(sv)^2(yz)^k-(tw)^2(zx)^k,
\end{align*}
because
\[
u^2-4ty=v^2,\qquad u^2-4sx=w^2.
\]
Similarly, the second formula in \eqref{eq:PRtransform} gives
\begin{align*}
P
&=r^2u^4(xy)^k
-s^2(u^2-4ty)^2(yz)^k
-t^2(u^2-4sx)^2(zx)^k\\
&=r^2u^4(xy)^k-s^2v^4(yz)^k-t^2w^4(zx)^k.
\end{align*}
Thus the transformed matrix has exactly the entries $P,Q,R$ introduced in Section~1.

\begin{theorem}[Fillipi identity]\label{thm:fillipi}
For every integer $n\ge2$,
\[
\mathcal A^2+\mathcal B^2-\mathcal C^2
=
(8rst)^2(xyz)^{n-2}(x^n+y^n-z^n).
\]
\end{theorem}

\begin{proof}
By \eqref{eq:baseidentity},
\[
\det\Gamma_0=-r^2(xyz)^kF.
\]
Since $\det T=-4st$, equation \eqref{eq:congruence} gives
\[
\det
\begin{pmatrix}
P&R\\
R&Q
\end{pmatrix}
=
16s^2t^2\det\Gamma_0
=
-16r^2s^2t^2(xyz)^kF.
\]
Finally,
\[
\mathcal A^2+\mathcal B^2-\mathcal C^2
=
-4(PQ-R^2),
\]
which proves the assertion.
\end{proof}

In particular, every solution of $x^n+y^n=z^n$ satisfies
\[
\mathcal A^2+\mathcal B^2=\mathcal C^2.
\]

\begin{remark}\label{rem:quartic}
The occurrence of the particular fourth powers $u^4,v^4,w^4$ is not intrinsic to the determinant identity. The simpler pair $(e_3,e)$ already gives \eqref{eq:baseidentity}; the factor $16s^2t^2$ and the fourth-power terms in Fillipi's expressions arise from the polynomial column transformation \eqref{eq:basis}.
\end{remark}

\section{Rank-one factorisation on the Fermat hypersurface}

The same matrix formulation determines the exact algebraic content of the Pythagorean equation obtained from Theorem~\ref{thm:fillipi}. Put
\begin{equation}\label{eq:pq}
q=x^{n-1}+y^{n-1}+z^{n-1},
\qquad
p=u^2q-4stz^{n-1}.
\end{equation}

\begin{proposition}\label{prop:rankone}
Modulo $F=x^n+y^n-z^n$,
\begin{equation}\label{eq:rankone}
\begin{pmatrix}
P&R\\
R&Q
\end{pmatrix}
\equiv
-
\begin{pmatrix}
p\\ q
\end{pmatrix}
\begin{pmatrix}
p&q
\end{pmatrix}.
\end{equation}
Consequently, at every point of $F=0$ over a commutative ring,
\begin{equation}\label{eq:pqr}
P=-p^2,\qquad Q=-q^2,\qquad R=-pq,
\end{equation}
and hence
\begin{equation}\label{eq:euclid}
\mathcal A=q^2-p^2,\qquad
\mathcal B=-2pq,\qquad
\mathcal C=-(p^2+q^2).
\end{equation}
\end{proposition}

\begin{proof}
For $B=[g,e]$, equation \eqref{eq:matrixform} reads
\[
\Gamma_D(\xi;B)
=
F B^{\mathsf T}DB
-
(B^{\mathsf T}D\xi)(B^{\mathsf T}D\xi)^{\mathsf T}.
\]
Modulo $F$, only the second term remains. Further,
\[
\ip{\xi}{e}_D
=
x^{n-1}+y^{n-1}+z^{n-1}
=q
\]
and, by \eqref{eq:g},
\[
\ip{\xi}{g}_D
=
u^2q-4stz^{n-1}
=p.
\]
This proves \eqref{eq:rankone}; equations \eqref{eq:pqr} and \eqref{eq:euclid} follow.
\end{proof}

\begin{remark}
The same column transformation also explains the parameter $p$. For the base pair $B_0=[e_3,e]$,
\[
B_0^{\mathsf T}D\xi
=
\begin{pmatrix}
z^{n-1}\\ q
\end{pmatrix}.
\]
Since $B=B_0T$,
\[
B^{\mathsf T}D\xi
=
T^{\mathsf T}
\begin{pmatrix}
z^{n-1}\\ q
\end{pmatrix}
=
\begin{pmatrix}
u^2q-4stz^{n-1}\\ q
\end{pmatrix}.
\]
Thus both the multiplicative factor in Theorem~\ref{thm:fillipi} and the polynomial $p$ in \eqref{eq:pq} arise from the same $2\times2$ column transformation.
\end{remark}

\begin{remark}\label{rem:arithmetic}
Equation \eqref{eq:euclid} is the standard Euclidean form of a Pythagorean triple, up to signs. Since it holds identically in the quotient ring
\[
\mathbb Z[x,y,z]/(x^n+y^n-z^n),
\]
the associated Pythagorean system is a universal polynomial consequence of the Fermat-type equation. Its solvability therefore supplies no additional algebraic condition on $x,y,z$; no converse is asserted.
\end{remark}

\end{document}